\documentclass[preprint,12pt,3p]{elsarticle}

\usepackage{amssymb}
\usepackage{amsthm}

\usepackage{amsmath}
\usepackage{mathtools}
\usepackage{amsfonts}
\usepackage{mathrsfs}
\usepackage{verbatim}
\usepackage{etoolbox}
\usepackage[none]{hyphenat}

\newtheorem{thm}{Theorem}

\newtheorem{rmk}[thm]{Remark}

\journal{arXiv}

\begin{document}

\begin{frontmatter}

\title{An Identity for the Partition Function Involving Parts of $k$ Different Magnitudes\tnoteref{label1}}
\tnotetext[label1]{This research did not receive any specific grant from funding agencies in the public, commercial, or not-for-profit sectors.}

\author{Saud Hussein}
\address{Institute of Mathematics, Academia Sinica, 6F, Astronomy-Mathematics Building, No.1, Sec.4, Roosevelt Road, Taipei 10617, Taiwan}

\ead{saudhussein@gate.sinica.edu.tw}

\begin{abstract}
Using previous work by Merca, we show the partition function involving parts of $k$ different magnitudes, shifted by the triangular numbers $\binom{k+1}{2}$, equals the self convolution of the unrestricted partition function. We also provide a combinatorial proof of this result.
\end{abstract}


\end{frontmatter}


\section{Introduction}

Recall that a \textit{partition} of a positive integer $n$ is a non-increasing sequence of positive integers whose sum is $n$. The partition function $p(n)$ counts the number of such partitions of $n$. For example, since the partitions of $4$ are \[4, \quad 3 + 1, \quad 2 + 2, \quad 2 + 1 + 1, \quad \text{and} \quad 1+1+1+1,\] $p(4) = 5$. We also take $p(0) = 1$.

Euler \cite{Euler} began the mathematical theory of partitions in 1748 and MacMahon \cite{MacMahon} in 1921 was the first to study the number of partitions of $n$ that have exactly $k$ different values for the parts. Let $p(k,n)$ denote this function. For example, $p(3,8) = 5$ since the five partitions in question are \[5+2+1, \quad 4+3+1, \quad 4+2+1+1, \quad 3+2+2+1, \quad 3+2+1+1+1.\] Let $q(k,n)$ denote the number of partitions of $n$ into exactly $k$ distinct parts. So $q(3,8) = 2$ since we have the two partitions \[5+2+1, \quad 4+3+1.\]
Clearly $q(k,n) \leq p(k,n)$ and $p(k,n) = 0$ when $n < \binom{k+1}{2}$. We also take $q(0,0) = p(0,0) = 1$ and $p(0,n) = p(k,0) = 0$ for $k, n > 0$.

Due to Merca \cite[Corollary 1.2]{Merca2}, we have the following: Let $k,n \geq 0$ be integers. Then \begin{align} p(k,n) = \sum_{m = \binom{k+1}{2}}^n a(k,m) \, p(n-m) \label{merc1} \end{align} where \begin{align} a(k,m) = \sum_{j=k}^\infty (-1)^{j-k}\binom{j}{k} q(j,m). \label{merc2} \end{align}
In this paper, we use this result to prove the following:

\begin{thm} \label{MAIN}
Let $k,n \geq 0$ be integers with $k \geq n$. Then \[p\left(k, n + \binom{k+1}{2} \right) = \sum_{m = 0}^n p(m) \, p(n-m).\]
\end{thm}

\begin{rmk}
The sequence of integers, $(A_n)_{n\geq 0}$, given by the self convolution \[A_n =  \sum_{m = 0}^n p(m) \, p(n-m),\] is listed in \cite[sequence number A000712]{OEIS} and may be thought of as the number of partitions of $n$ into parts of two kinds. For example, $A_3 = 10$ since we have \[3, \quad \overline{3}, \quad 2+1, \quad \overline{2}+\overline{1}, \quad \overline{2}+1, \quad 2+\overline{1},\]
\[1+1+1, \quad \overline{1}+\overline{1}+\overline{1}, \quad \overline{1}+1+1, \quad 1+\overline{1}+\overline{1}.\]

Let $q(n)$ denote the number of partitions of $n$ into distinct parts. Then by Merca \cite[Theorem 5.2]{Merca3}, \[\sum_{k=0}^n p\left(n- \binom{k+1}{2}\right) = \sum_{m=0}^n q(m) \, q(n-m),\] a result much like Theorem~\ref{MAIN}.
Also, from Merca \cite[equation 12]{Merca1}, we have \begin{align} q(k,n) = p_k\left(n - \binom{k+1}{2}\right), \label{merc3} \end{align} where $p_k(n)$ denotes the number of partitions of $n$ with no part greater than $k$. Clearly \[q\left(k,n + \binom{k+1}{2}\right) = p_k(n) = p(n)\] for $k \geq n$. Again, this identity is similar to Theorem~\ref{MAIN}.
\end{rmk}

\section{Proof of Theorem~\ref{MAIN}}

Let $k,n \geq 0$ be integers. Applying a change of variables to \eqref{merc1}, \begin{align} p\left(k,n+\binom{k+1}{2}\right) = \sum_{m = \binom{k+1}{2}}^{n + \binom{k+1}{2}} a(k,m) \, p\left(n + \binom{k+1}{2} - m\right). \label{merc4} \end{align}
Since $q(k+1,m) = 0$ for $m < \binom{k+2}{2} = \binom{k+1}{2} + k+1$, then by \eqref{merc2} and \eqref{merc3}, \[a(k,m) = q(k,m) = p_k\left(m - \binom{k+1}{2}\right) = p\left(m - \binom{k+1}{2}\right)\] for $m \leq \binom{k+1}{2} + k$. Thus by \eqref{merc4}, \begin{align*} p\left(k,n+\binom{k+1}{2}\right) &= \sum_{m = \binom{k+1}{2}}^{n + \binom{k+1}{2}} p\left(m - \binom{k+1}{2}\right) p\left(n + \binom{k+1}{2} - m\right)\\
&= \sum_{m = 0}^n p(m) \, p(n-m) \end{align*} for $k \geq n$ and the proof is complete.

\section{Combinatorial proof of Theorem~\ref{MAIN}}

Let $k,n \geq 0$ be integers with $n \leq k$. For $m \in [0,n]$, take a partition $a_1+\cdots +a_r$ of $m$ with $1\leq a_1 \leq \cdots \leq a_r$ and a partition $b_1+\cdots +b_s$ of $n-m$ with $1\leq b_1 \leq \cdots \leq b_s$. Then \[r+s \leq m+s \leq m+(n-m) = n\leq k,\] so $m \leq k-s$ and \begin{align*} n+\binom{k+1}{2} &= m +(n-m) + 1+2+\cdots+k\\
&=a_1 + \cdots + a_r + b_1 + \cdots + b_s + 1 + 2 +\cdots + k\\
&=a_1 + \cdots + a_r + c_1 + \cdots + c_k
\end{align*}
where \[c_j = \begin{cases} 
      j, & 1\leq j \leq k-s \\
      j+b_{j-(k-s)}, & k-s+1\leq j\leq k 
   \end{cases}.\]
Clearly, the parts $a_1,\dots,a_r$ are contained in $\{c_1,c_2,\dots,c_{k-s}\}$ and $c_1 < c_2 < \cdots < c_k$, so each pair of partitions of $m$ and $n-m$ uniquely defines a partition \[a_1 + \cdots + a_r + c_1 + \cdots + c_k\] of $n+ \binom{k+1}{2}$ having exactly $k$ different values for the parts. The product $p(m)p(n-m)$ then counts the number of partitions that may be constructed in this fashion for each $m \in [0,n]$, thus completing the proof.

The idea of this proof was suggested by a referee reading a previous version of this paper. We also see that setting $c_j = \overline{b_j}$ for $1\leq j \leq s$ would instead give us the partition \[a_1 + \cdots + a_r + \overline{b_1} + \cdots + \overline{b_s}\] of $n$ into parts of two kinds. This justifies the earlier remark on sequence $(A_n)_{n\geq 0}$ and provides some motivation for the given combinatorial proof.

The combinatorial proof gives us a simple way of writing down the partitions of $n+ \binom{k+1}{2}$ having exactly $k$ different values for the parts for any $k \geq n$. All we require are the partitions of $m$ for each $m \in [0,n]$. For example, Theorem~\ref{MAIN} implies there are ten partitions of $3+ \binom{k+1}{2}$ having exactly $k$ different values for the parts for any $k \geq 3$. The case $k=3$ is shown below.

\begin{center}
\begin{tabular}{ |c||c|c|c|c|c|c| } 
 \hline
 $m$ & $r$ & $s$ & $a_1 + \cdots + a_r$ & $b_1 + \cdots + b_s$ & $c_1 + c_2 + c_3$ & partitions\\ 
 \hline
 \hline
 0 & 0 & 1 & 0 & 3 & $1+2+6$ & $1+2+6$\\ [0.5ex]
 & 0 & 2 & 0 & $1+2$ & $1+3+5$ & $1+3+5$\\ [0.5ex]
 & 0 & 3 & 0 & $1+1+1$ & $2+3+4$ & $2+3+4$\\ [0.5ex]
 \hline
 1 & 1 & 1 & 1 & 2 & $1+2+5$ & $1+1+2+5$ \\ [0.5ex]
 & 1 & 2 & 1 & $1+1$ & $1+3+4$ & $1+1+3+4$\\ [0.5ex]
 \hline
 2 & 1 & 1 & 2 & 1 & $1+2+4$ & $1+2+2+4$\\ [0.5ex]
 & 2 & 1 & $1+1$ & 1 & $1+2+4$ & $1+1+1+2+4$\\ [0.5ex]
 \hline
 3 & 1 & 0 & 3 & 0 & $1+2+3$ & $1+2+3+3$\\ [0.5ex]
 & 2 & 0 & $1+2$ & 0 & $1+2+3$ & $1+1+2+2+3$\\ [0.5ex]
 & 3 & 0 & $1+1+1$ & 0 & $1+2+3$ & $1+1+1+1+2+3$\\ [0.5ex]
 \hline
\end{tabular}
\end{center}

\section*{Acknowledgment} The identity in this paper was discovered using Mathematica 11 and the author thanks the Institute of Mathematics, Academia Sinica for providing this tool.

\bibliographystyle{elsarticle-harv}

\bibliography{biblio.bib}

\end{document}